\documentclass{article}
\usepackage{amsmath}
\usepackage{amssymb}

\title{On the classification of simple modules for cyclotomic 
Hecke algebras of type $G(m\!,1\!,n)$ \\
and \\
Kleshchev multipartitions}
\author{Susumu Ariki}

\date{}

\newtheorem{thm}{Theorem}[section]
\newtheorem{defn}[thm]{Definition}

\newtheorem{lemma}[thm]{Lemma}

\newcommand{\Hn}{${\mathcal H_n}$}
\newcommand{\KP}{\hphantom{}^\gamma{\mathcal KP}}
\begin{document}
\maketitle

\begin{abstract}
We give a proof of a conjecture that 
Kleshchev multipartitions are those partitions which 
parametrize non-zero simple modules obtained as 
factor modules of Specht modules by their own radicals. 
\end{abstract}

\footnotetext
{A.M.S. subject classification, 20C20, 20C33, 20G05 
\\
This work is a contribution to the JSPS-DFG Japanese-German 
Cooperative Science Promotion Program on 
"Representation Theory of Finite and Algebraic Groups"
}

\section{Introduction}

After Hecke algebras appeared, unexpectedly deep 
applications and results have been found in 
the representation theory of these algebras. 
Concerned with ordinary representations, Lusztig's cell theory 
is the main driving force. But we do not 
consider it here. The other interest is about the 
modular representation theory of these algebras. 
We are mainly working with Hecke algebras of type 
$A$ and type $B$, and this research is driven by 
Dipper and James \cite{DJ1}\cite{DJ2}.  
Recently, a new type of Hecke algebras was introduced. 
We call them cyclotomic Hecke algebras 
of type $G(m\!,1\!,n)$ following \cite{BM}. 
Hecke algebras of type $A$ and type $B$ are 
special cases of these algebras. The author studied 
modular representations over the algebra for 
the case that parameters were 
roots of unity in the field of complex numbers \cite{A1}. 
In particular, it gives the 
classification of simple modules. The removal of the restriction 
on base fields was achieved in \cite{AM}. 
In the paper \cite{AM}, we gave a 
classification of the simple modules over cyclotomic 
Hecke algebras in terms of 
the crystal graphs of integrable highest weight modules 
over certain quantum algebras. The result turns out to be useful for 
verifying a conjecture of Vigneras \cite{Vig3}. 

On the other hand, another approach was already proposed 
in \cite{GL}\cite{DJM}. Main results in the theory are that 
we can define "Specht modules", and that 
each Specht module $S^{\underline\lambda}$ has natural bilinear form, and 
each of $D^{\underline\lambda}:=S^{\underline\lambda}
/{\rm rad}S^{\underline\lambda}$ is an absolutely 
irreducible or zero module. Further, the theory claims that 
the set of non-zero $D^{\underline\lambda}$ is a 
complete set of simple modules. 

But there is one drawback. The theory does not tell 
which $D^{\underline\lambda}$ are actually non-zero. 
We conjectured in \cite{AM} that the crystal graph description gave the 
criterion. Namely, 
we conjectured that $D^{\underline\lambda}\ne0$ if and only if 
$\underline\lambda$ is a Kleshchev multipartition. The purpose of this 
paper is to prove the conjecture. It is achieved by interpreting 
the conjecture into a problem about canonical bases in Fock spaces. 
This part is based on \cite{A1} and \cite{AM}. 
Then it is easily verified by using a recent result of Uglov \cite{U}. 

The author is grateful to A.Mathas for 
discussion he had at the early stage of the research. 
He also thanks B.Leclerc, Varagnolo and Vasserot. 

\section{Preliminaries}

Let $R$ be an integral domain, $u_1,\dots,u_m$ be elements 
in $R$, and $\zeta$ be an invertible element. 
The Hecke algebra of type $G(m,1,n)$ is the algebra 
associated with these parameters is the $R$-algebra 
defined by the following defining relations for generators 
$a_i$ $(1\!\le\! i\!\le\! n)$. We denote this algebra by 
\Hn. 

\[
\begin{array}{c}
(a_1-u_1)\cdots(a_1-u_m)=0,\qquad (a_i-\zeta)(a_i+\zeta^{-1})=0\quad(i\ge 2)\\
\\
a_1a_2a_1a_2=a_2a_1a_2a_1, \qquad a_ia_j=a_ja_i\quad(j\ge i\!+\!2)\\
\\
a_ia_{i-1}a_i=a_{i-1}a_ia_{i-1}\;(3\!\le\! i\!\le\! n)
\end{array}
\]

It is known that this algebra is $R$-free of rank $m^nn!$ as an $R$-module. 
This algebra is also known to be cellular in the sense of 
Graham and Lehrer \cite{GL}, and thus has Specht modules. 
Following \cite{DJM}, we shall explain the theory. 
A {\bf partiton} $\lambda$ of size $n$ is a sequence of non-negative integers 
$\lambda_1\ge\lambda_2\ge\cdots$ such that $\sum \lambda_i=n$. 
We write $|\lambda|=n$. A {\bf multipartiton} of size $n$ is a 
sequence of $m$ partitions 
$\underline\lambda=(\lambda^{(1)},\dots,\lambda^{(m)})$ 
such that $\sum_{k=1}^m|\lambda^{(k)}|=n$. If $n=0$, we denote 
the multipartition by $\underline\emptyset$. 
The set of multipartitions has a poset structure. 
The partial order is the {\bf dominance order}, which is 
defined as follows. 

\begin{defn}
Let $\underline\lambda$ and $\underline\mu$ be multipartitions. 
We say that $\underline\lambda$ dominates $\underline\mu$, and 
write $\underline{\lambda}\trianglerighteq\underline{\mu}$ if we have 
\[
\sum_{l=1}^{k-1} |\lambda^{(l)}|+\sum_{i=1}^j \lambda^{(k)}_i
\ge 
\sum_{l=1}^{k-1} |\mu^{(l)}|+\sum_{i=1}^j \mu^{(k)}_i
\]
for all $j, k$. 
\end{defn}

With each multipartition 
$\underline\lambda$, we can associate an \Hn-module 
$S^{\underline\lambda}$. Its concrete construction is explained in 
\cite[(3.28)]{DJM}. It is easy to see from the construction 
that it is free as an $R$-module. 
These modules are 
called {\bf Specht modules}. 
Each Specht module is naturally equipped with 
a bilinear form \cite[(3.28)]{DJM}. 
We set $D^{\underline\lambda}=S^{\underline\lambda}/
{\rm rad}\,S^{\underline\lambda}$. It can be zero, 
but non-zero ones exhaust all simple \Hn-modules. 
We denote the projective cover of $D^{\underline\lambda}$ by 
$P^{\underline\lambda}$. 

We remark that Graham and Lehrer have introduced the notion 
of cellular algebras and have developped general theory for 
classifying simple modules by using "cell modules". 
In \cite{GL}, the cellular bases for the cell modules 
are given by Kazhdan-Lusztig bases. Here, 
different cellular bases are given, but the strategy to classify 
simple modules is the same. Hence we call the following 
parametrization the Graham-Lehrer parametrization. 

\begin{thm}[{\cite[Theorem 3.30]{DJM}}]
\label{GL parametrization} 
Suppose that $R$ is a field. Then, 

\begin{description}
\item[{\rm (1)}]
Non-zero $D^{\underline\lambda}$ form 
a complete set of non-isomorphic simple \Hn-modules. Further, 
these modules are absolutely irreducible. 
\item[{\rm (2)}]
Let $\underline\lambda$ and $\underline\mu$ be multipartitions of $n$ and 
suppose that $D^{\underline\mu}\ne 0$ and that 
$[S^{\underline\lambda}:D^{\underline\mu}]\ne 0$. 
Then we have $\underline{\lambda}\trianglerighteq\underline{\mu}$. 
\item[{\rm (3)}]
$[S^{\underline\lambda}:D^{\underline\lambda}]=1$. 
\end{description}
\end{thm}

\noindent
Note that (2) is equivalent to the following (2').

\begin{description}
\item[{\rm\small (2')}]
{\it Let $\underline\lambda$ and $\underline\mu$ be multipartitions of $n$ and 
suppose that $D^{\underline\mu}\ne 0$ and that 
$[P^{\underline\mu}:S^{\underline\lambda}]\ne 0$. 
Then we have $\underline{\lambda}\trianglerighteq\underline{\mu}$. }
\end{description}

\noindent
It is obvious since we have 
$[P^{\underline\mu}:S^{\underline\lambda}]=
{\rm dim}{\it Hom}_{{\mathcal H_n}}
(P^{\underline\mu}, S^{\underline\lambda})=
[S^{\underline\lambda}:D^{\underline\mu}]$. 

As is explained in \cite[1.2]{AM}, the classification of 
simple \Hn-modules is reduced to the classification 
in the case that $u_1,\dots,u_m$ are powers of $\zeta^2$. 
This is a consequence of a result in \cite[2.13]{Vig1} 
(see also \cite{DM}). We can also assume that $\zeta^2\ne1$, since 
the case $\zeta^2=1$ is well understood. In the rest of the paper 
throughout, we assume that 
\[
u_i=\zeta^{2\gamma_i} \quad (\,i=1,\dots,m\,),\qquad \zeta^2\ne1
\]
If $\zeta^2$ is a primitive 
$r$ th root of unity for a natural number $r$, 
$\gamma_i$ take values in ${\mathbb Z}/r{\mathbb Z}$. 
Otherwise, these take values in ${\mathbb Z}$. 
We now recall the notion of{\bf Kleshchev multipartitions} 
associated with $(\gamma_1,\dots,\gamma_m)$. To do this, 
we explain the notion of good nodes first. 

We identify a multipartition $\underline\lambda=
(\lambda^{(1)},\dots,\lambda^{(m)})$ with 
the associated Young diagram, i.e. an $m$-tuple of 
the Young diagrams 
associated with $\lambda^{(1)},\dots,\lambda^{(m)}$. 
Let $x$ be a node on the Young diagram which is 
located on the $a$ th row and the $b$ th column of 
$\lambda^{(c)}$. If $u_c\zeta^{2(b-a)}=
\zeta^{2i}$, we say that the node $x$ has 
{\bf residue} $i$ (with respect to 
$\gamma=(\gamma_1,\dots,\gamma_m)$). 
We denote the residue by $r_\gamma(x)$. A node is 
called an {\bf $i$-node} 
if its residue is $i$. Let $\underline\lambda$ and 
$\underline\mu$ be multipartitions. We first assume that 
$|\underline\lambda|\!+\!1=|\underline\mu|$, $r_\gamma(x)\equiv i$, 
and let $x$ be 
the node $\underline\mu/\underline\lambda$. We then call 
$x$ an {\bf addable $i$-node} of $\underline\lambda$. 
If $|\underline\lambda|\!-\!1=
|\underline\mu|$ and $r_\gamma(x)\equiv i$, we call 
$x=\underline\lambda/\underline\mu$ 
a {\bf removable $i$-node} of $\underline\lambda$. 

For each residue $i$, we have the notion of normal $i$-nodes 
and good $i$-nodes. To define these, 
We read addable 
and removable $i$-nodes of $\underline\lambda$ in the 
following way. We start with 
the first row of $\lambda^{(1)}$, and we read rows in 
$\lambda^{(1)}$ downward. We then move to the first row 
of $\lambda^{(2)}$, and repeat the same procedure. We continue 
the procedure to $\lambda^{(3)},\dots,\lambda^{(m)}$. 
If we write $A$ for an addable $i$-node, and similarly $R$ for a 
removable $i$-node, we get a sequence of $A$ and $R$. We then 
delete $RA$ as many as possible. For example, if the sequence is 
$RRAAAARRRAARAR$, it ends up with 
$----AAR------R$. The remaining removable $i$-nodes in this sequence 
are called {\bf normal $i$-nodes}. The node corresponding 
to the leftmost $R$ is called the {\bf good $i$-node}. If 
$x$ is a good $i$-node for some $i$, we simply say $x$ is a 
{\bf good node}. 
We can now define the set of Kleshchev multipartitions associated 
with $\gamma=(\gamma_1,\dots,\gamma_m)$. 

\begin{defn}
We declare that $\underline\emptyset$ is Kleshchev. 
Assume that we have already defined 
the set of Kleshchev multipartitions of size $n$. 

Let $\underline\lambda$ be a multipartition of size $n\!+\!1$. 
We say that $\underline\lambda$ is Kleshchev 
if and only if there is a good node $x$ of $\underline\lambda$ 
such that $\underline\mu:=\underline\lambda\setminus\{x\}$ is a Kleshchev 
multipartition. 
\end{defn}

We denote the set of Kleshchev multipartitions of 
size $n$ by $\KP_n$, and set $\KP=\sqcup_{n\ge0}\KP_n$. 

\begin{thm}[{\cite[Theorem C]{AM}}]
\label{AM parametrization}
Suppose that $\zeta^2$ and $u_i$ satisfy the above 
condition. 
Then, the irreducible \Hn-modules are indexed by the set of 
Kleshchev multipartitions. 
\end{thm}

Hence we have two parametrizations. 
One given in Theorem \ref{AM parametrization} and 
one given in Theorem \ref{GL parametrization}. 
It is natural to ask, if these 
coincide. The main observation is the 
following conjecture, 
which will be proved in the 
last section. The conjecture was formulated by Mathas. 

\bigskip
\noindent
{\bf Conjecture}\cite[2.12]{AM}\;
These two parametrizations coincide. In particular, 
$D^{\underline\lambda}\ne0$ if and only if 
$\underline\lambda$ is a Kleshchev multipartition. 

\bigskip
To prove this, we use certain Fock spaces over a quantum 
algebra. \footnote{The idea to use such Fock spaces to study 
the modular representation theory of 
cyclotomic Hecke algebras first appeared in \cite{A1}, 
generalizing and verifying a conjecture of Lascoux, Leclerc 
and Thibon \cite{LLT}.} In the next section, we recall necessary 
ingredients of these Fock spaces. 

\section{Fock space}

Recall that the multiplicative order of $\zeta^2$ is $r\ge2$. 
We denote by $U_v$ the quantum algebra of type $A_{r-1}^{(1)}$ if 
$r$ is finite, and of type $A_\infty$ if $r\!=\!\infty$. 
Let ${\mathcal F}_v^\gamma$ be the combinatorial Fock space:
it is a $U_v$-module, whose basis elements are 
indexed by the set of all multipartitions. 
We identify the basis elements 
with the multipartitions. The size of multipartitions 
naturally makes it into a graded vector space. 

We consider the $U_v$-submodule ${\mathcal M}_v^\gamma$ of 
${\mathcal F}_v^\gamma$ 
generated by $\underline\emptyset$. It is isomorphic to 
an irreducible highest weight module with highest weight 
$\Lambda=\Lambda_{\gamma_1}\!+\!\cdots\!+\!\Lambda_{\gamma_m}$, 
where $\Lambda_i$ are fundamental weights. 
To describe its basis in a combinatorial way, we need 
the crystal graph theory of Kashiwara. In our particular setting, 
we can prove the following theorem using 
argument in \cite{MM}. The theorem explains 
the representation theoretic meaning of Kleshchev multipartitions. 

\begin{thm}[{\cite[Theorem 2.9, Corollary 2.11]{AM}}]
\label{Misra Miwa}
Let $R_v$ be the localized ring of ${\mathbb Q}[v]$ with respect to 
the prime ideal $(v)$. We consider the $R_v$-lattice of 
${\mathcal F}_v^\gamma$ generated by all multipartitions, and denote 
it by ${\mathcal L}_v^\gamma$. 
We set $L(\Lambda)=
{\mathcal L}_v^\gamma\cap {\mathcal M}_v^\gamma$, and 
$B(\Lambda)=\{ \underline\lambda\;mod\;v\,|\,\underline\lambda\in 
\KP\,\}$. Then, $(L(\Lambda),B(\Lambda))$ is a (lower) 
crystal base of ${\mathcal M}_v^\gamma$ 
in the sense of Kashiwara. 
\end{thm}

To explain the $U_v$-module structure 
given to ${\mathcal F}^\gamma$, we first fix 
notations. 
Let $\underline\lambda$ be a multipartition and 
let $x$ be a node on the associated Young diagram 
which is located on the $a$ th row 
and $b$ th column of $\lambda^{(c)}$. Then 
we say that a node is {\bf above} $x$ if it is on $\lambda^{(k)}$ 
for some $k<c$, or if it is on $\lambda^{(c)}$ and the row number is 
strictly smaller than $a$. We denote the set of 
addable (resp. removable) $i$-nodes of $\underline\lambda$ 
which are above 
$x$ by $A_i^a(x)$ (resp. $R_i^a(x)$). In a similar way, we say that 
a node is {\bf below} $x$ if it is on $\lambda^{(k)}$ 
for some $k>c$, or if it is on $\lambda^{(c)}$ and the row number is 
strictly greater than $a$. We denote 
the set of addable (resp. removable) $i$-nodes of $\underline\lambda$ 
which are 
below $x$ by $A_i^b(x)$ (resp. $R_i^b(x)$). 
The set of all addable (resp. removable) $i$-nodes of 
$\underline\lambda$ is 
denoted by $A_i(\underline\lambda)$ (resp. $R_i(\underline\lambda)$). 

In the similar way, we define the notion that a node is 
{\bf left} to $x$ (resp. {\bf right} to $x$). We denote 
the set of addable $i$-nodes which are left to $x$ 
(resp. right to $x$) by 
$A_i^l(x)$ (resp. $A_i^r(x)$). The set of removable $i$-nodes 
which are left to $x$ (resp. right to $x$) is denoted by 
$R_i^l(x)$ (resp. $R_i^r(x)$). 
We then set 
\[
\begin{array}{c}
N_i^a(x)=|A_i^a(x)|-|R_i^a(x)|,\quad
N_i^b(x)=|A_i^b(x)|-|R_i^b(x)|\\
\\
N_i(\underline\lambda)=|A_i(\underline\lambda)|-
|R_i(\underline\lambda)|
\end{array}
\]
$N_i^l(x)$ and $N_i^r(x)$ are similarly defined. 
Finally, we denote the number of all $0$-nodes in 
$\underline\lambda$ by $N_d(\underline\lambda)$. Then 
the $U_v$-module structure of ${\mathcal F}_v^\gamma$ 
(called Hayashi action) is defined as follows. 
\[
\begin{array}{c}
e_i{\underline\lambda}=\displaystyle\sum_{r_\gamma
(\underline\lambda/\underline\mu)\equiv i}
v^{-N_i^a(\underline\lambda/\underline\mu)}\underline\mu, \quad
f_i{\underline\lambda}=\displaystyle\sum_{r_\gamma
(\underline\mu/\underline\lambda)\equiv i}
v^{N_i^b(\underline\mu/\underline\lambda)}\underline\mu \\
\\
v^{h_i}\underline\lambda=v^{N_i(\underline\lambda)}\underline\lambda, \quad
v^d\underline\lambda=v^{-N_d(\underline\lambda)}\underline\lambda
\end{array}
\]
Let ${\mathcal F}_v^{\gamma_i}$ $(i=1,\dots,m)$ be the 
combinatorial Fock spaces defined 
for the cases that $m\!=\!1$ and $\gamma=\gamma_i$. 
The $U_v$-module structure on these spaces are defined by 
the same formula given above. 
Let $\Delta'$ be the comultiplication defined by 
\[
\Delta'(e_i)=v^{-h_i}\otimes e_i + e_i\otimes1, \quad
\Delta'(f_i)=1\otimes f_i + f_i\otimes v^{h_i}
\]
If we identify ${\mathcal F}_v^\gamma$ with 
${\mathcal F}_v^{\gamma_1}\otimes\cdots\otimes
{\mathcal F}_v^{\gamma_m}$, the representation on ${\mathcal F}_v^\gamma$ 
coincides with the tensor product representation defined 
with respect to $\Delta'$. 

Let 
$-\gamma':=(-\gamma_m,\dots,-\gamma_1)$. Namely, we write 
$r_{-\gamma'}(x)\equiv i$ if 
a node $x$ on the $a$ th row and the $b$ th column of 
the $\lambda^{(c)}$ satisfies 
$u_{m\!+\!1\!-\!c}^{-1}\zeta^{2(b-a)}=\zeta^{2i}$. 

For each partition $\lambda$ we denote its transpose by 
$\lambda'$. For a multipartition $\underline\lambda$, 
we denote $({\lambda^{(m)}}',\dots,{\lambda^{(1)}}')$ by 
$\underline\lambda'$ and call it the {\bf flip transpose} of 
$\underline\lambda$. Similarly, we denote 
$({\lambda^{(1)}}',\dots,{\lambda^{(m)}}')$ by $\underline\lambda^{\rm T}$ 
and call it the {\bf transpose} of $\underline\lambda$. 

Let $\sigma:{\mathcal F}_v^{-\gamma'}\rightarrow{\mathcal F}_v^\gamma$ 
be a linear map 
which maps $\underline\lambda$ to 
$\underline\lambda'$. Then the coproduct on 
${\cal F}_v^{-\gamma'}$ coincides with Kashiwara's, 
and the action coincides with \cite{LLT}. Hence the $R_v$-lattice 
generated by $\underline\lambda$ is a crystal lattice in 
${\mathcal F}_v^\gamma$. 

Let $\xi:{\mathcal F}_v^{\gamma}\rightarrow
{\mathcal F}_{v^{-1}}^{-\gamma}$ 
be a semilinear map which sends $\underline\lambda$ to 
$\underline\lambda^{\rm T}$. Then we have a representation 
which is compatible with Lusztig's coproduct. 
The space ${\mathcal F}_{v^{-1}}^{-\gamma}$ is the same space as 
${\mathcal F}_v^{-\gamma}$, but to stress that the crystal base 
here is a so-called "basis at $v=\infty$" in the sense of Lusztig, 
and not the one generated over $R_v$ by $\underline\lambda$, 
we adopt the different notation. 
The action on ${\mathcal F}_{v^{-1}}^{-\gamma}$ is as follows. 
We also call it Hayashi action. 
\[
\begin{array}{c}
e_i{\underline\lambda}=\displaystyle\sum_{r_\gamma
(\underline\lambda/\underline\mu)\equiv i}
v^{N_i^l(\underline\lambda/\underline\mu)}\underline\mu, \quad
f_i{\underline\lambda}=\displaystyle\sum_{r_\gamma
(\underline\mu/\underline\lambda)\equiv i}
v^{-N_i^r(\underline\mu/\underline\lambda)}\underline\mu \\
\\
v^{h_i}\underline\lambda=v^{N_i(\underline\lambda)}\underline\lambda, \quad
v^d\underline\lambda=v^{-N_d(\underline\lambda)}\underline\lambda
\end{array}
\]
In the rest of paper, we exclusively work with 
${\mathcal F}_{v^{-1}}^{-\gamma}$. 

\section{The proof of the conjecture}

We first interprete the conjecture into a problem about canonical 
bases on Fock spaces. To do this, we use 
the direct sum of the Grothendieck groups of projective 
\Hn-modules. 
We always assume that the coefficients are 
extended to the field of rational numbers. 
If \Hn is semisimple, all $S^{\underline\lambda}$ 
are irreducible, and we identify the direct sum with 
${\mathcal F}_{v\!=\!1}^\gamma$, which is by definition 
a based ${\mathbb Q}$-vector space 
whose basis elements are indexed by multipartitions, and nodes of 
multipartitions are given residues. If \Hn is 
not semisimple, we have a proper subspace of 
${\mathcal F}_{v\!=\!1}^\gamma$ by lifting idempotents argument. 
It is proved in \cite{A1} that 
it coincides with ${\mathcal M}_{v\!=\!1}^\gamma$. 

Recall that simple modules are 
obtained as factor modules of Specht modules. To distinguish 
between simple modules over different base rings, we write 
$D^{\underline\lambda}_R$ when the base ring is $R$. 
Let $(K,R,F)$ be a modular system. We assume that there is 
an invertible element $\zeta\in R$ such that its multiplicative order in 
$K$ and $F$ is the same. 
Then $D^{\underline\lambda}_K$ 
is obtained from $D^{\underline\lambda}_R$ 
by extension of coefficients, and 
$D^{\underline\lambda}_F$ is obtained from 
$D^{\underline\lambda}_R$ by taking the unique simple factor 
module of $D^{\underline\lambda}_R\otimes F$. 
The proof of Theorem \ref{AM parametrization} implies that 
these give the correspondence between simple modules over fields of 
positive characteristics and fields of characteristic $0$, and 
$D^{\underline\lambda}_F\ne 0$ 
if and only if $D^{\underline\lambda}_K\ne 0$. 
Further, still assuming that the multiplicative order is the same, 
the proof given in \cite{AM} also shows 
that $D^{\underline\lambda}_K\ne 0$ 
if and only if $D^{\underline\lambda}_{\mathbb C}\ne 0$. 
In particular, to know 
which $D^{\underline\lambda}$ are non-zero, 
it is enough to consider the case that the base field is ${\mathbb C}$. 

Now assume that we are in the case that the base field is 
${\mathbb C}$. We identify the direct sum of the Grothendieck 
groups of projective \Hn-modules with 
${\mathcal M}_{v\!=\!1}^\gamma$ as before. 
The main theorem in \cite{A1} asserts that 
the canonical basis evaluated at $v\!=\!1$ consists of 
indecomposable projective \Hn-modules $(n=0,1\dots)$. 
Hence we have a bijection between canonical basis elements of 
${\mathcal M}_v^\gamma$ and indecomposable projective \Hn-modules 
$P^{\underline\lambda}$ for various $n$, and thus a bijection between 
canonical basis elements of 
${\mathcal M}_v^\gamma$ and simple \Hn-modules 
$D^{\underline\lambda}$ for various $n$. 

It is known that the canonical basis 
gives a crystal base of ${\mathcal M}_v^\gamma$ \cite{G'L'}, 
which is unique up to scalar \cite{Ka}. More precisely, 
the crystal lattice $L(\Lambda)$ is the 
$R_v$-lattice generated by the canonical basis elements, and 
$B(\lambda)$ consists of the canonical basis elements modulo $v$. 
Then Theorem \ref{Misra Miwa} asserts that with each canonical 
basis element $G(b)$, we can uniquely associate a multipartition 
$\underline\nu \in\KP$. To summarize, we have the following. 

\begin{quotation}
\noindent
{\it For each non-zero $D^{\underline\lambda}$, there exists 
a unique canonical basis element $G(b)$ such that we have 
$G(b)_{v\!=\!1}=P^{\underline\lambda}$ and 
$G(b)\equiv \underline\nu \;mod\;vL(\Lambda)$.}
\end{quotation}

If $\underline\nu=\underline\lambda$ holds in general, then 
the following lemma proves the conjecture. 
The dominance order on ${\mathcal F}_{v^{-1}}^{-\gamma}$ is 
defined by reading columns from left to right. If 
we read the columns from right to left, we have the 
reversed dominance order. We 
denote it by $\underline\lambda\ge\underline\mu$.

\begin{lemma} 
\label{minimum} 
Assume that for every canonical basis element $G(b)$ in 
${\mathcal F}_v^\gamma$, there exists a unique maximal element 
among the multipartitions 
appearing in $\xi(G(b))$ with respect to the reversed 
dominance order, and assume that it has coefficient $1$. 
Then we have that 
$D^{\underline\lambda}\ne0$ if and only if 
$\underline\lambda$ is a Kleshchev multipartition. 
\end{lemma}
(Proof) Recall that ${\mathcal M}_{v\!=\!1}^\gamma$, the sum of 
Grothendieck groups of projective \Hn-modules, is embedded 
into ${\mathcal F}_{v\!=\!1}^\gamma$ by sending $S^{\underline\lambda}$ 
to $\underline\lambda$. Hence, Theorem \ref{GL parametrization} 
implies that $\xi(P^{\underline\lambda})$ has the form
$\xi(P^{\underline\lambda})=
\underline\lambda^{\rm T} + \sum_{\underline\mu<
\underline\lambda^{\rm T}} 
c_{\underline\mu}\,\underline\mu$. In particular, among 
multipartitons appearing in $\xi(P^{\underline\lambda})$, 
$\underline\lambda^{\rm T}$ is the maximal element 
with respect to the dominance order. 

On the other hand, Theorem \ref{AM parametrization} implies 
that there exists a canonical basis element $G(b)$ such that 
$P^{\underline\lambda}=G(b)_{v\!=\!1}$. We apply the assumption 
to $\xi(G(b))$. Then 
multipartitions appearing in 
$\xi(G(b))$ has a unique maximal element with 
coefficient $1$. 
Since it is a canonical basis element and the 
coefficient is $1$, it must be the 
transpose of a Kleshchev multipartiton, say $\underline\nu$. 
We specialize $\xi(G(b))$ to $v\!=\!1$. Note that $\underline\nu$ 
does not vanish. Since both $\underline\lambda^{\rm T}$ and 
$\underline\nu^{\rm T}$ are maximal elements, 
we have $\underline\nu=\underline\lambda$. 
Hence 
the two parametrizations given in Theorem \ref{GL parametrization} 
and Theorem \ref{AM parametrization} coincide. 
$\blacksquare$

\bigskip
We say that $\tilde{\gamma}=
(\tilde{\gamma}_1,\dots,\tilde{\gamma}_m)$ is a {\bf lift} 
of $\gamma=(\gamma_1,\dots,\gamma_m)$ if 
$\tilde{\gamma}_i\;mod\;r=\gamma_i$ for all $i$. If $r\!=\!\infty$, 
we set $\tilde{\gamma}=\gamma$. For each $\tilde{\gamma}$, 
we have ${\mathcal F}_{v^{-1}}^{-\tilde\gamma}$, and 
it is a module over the quantum algebra of type $A_\infty$. 
It will play a main role in the following. 

In \cite{TU}, Takemura and Uglov constructed higher level Fock spaces 
by generalizing \cite[Proposition 1.4]{KMS}. 
Let $\{\,u_i\,\}_{i\in{\mathbb Z}}$ be the basis vectors of 
an infinite dimensional space. More precisely, 
the space is originally ${\mathbb Q}(v)^r\otimes 
{\mathbb Q}(v)^m[z,z^{-1}]$, and if we denote 
the basis elements by $e_a\otimes e_bz^N$, we identify 
$u_i$ with $e_a\otimes e_bz^N$ through $i=a+r(b-1-mN)$ as in 
\cite{U}. Note that this identification is different from 
that in \cite{TU} since 
the evaluation representation for $U_v'(\hat{\mathfrak sl}_m)$ 
taken in \cite{TU} is different from that in \cite{U}. 
This space is naturally a $U_v'(\hat{\mathfrak sl}_r)\otimes 
U_v'(\hat{\mathfrak sl}_m)$-module. The semi-infinite wedges 
of the form 
$u_I=u_{i_1}\wedge u_{i_2}\wedge\cdots$ 
such that $i_k=c\!-\!k\!+\!1$ for all $k>\!>0$ are 
said to have charge $c$. The space 
of semi-infinite wedges of charge $c$ is denoted by 
${\cal F}_c$. To make ${\cal F}_c$ into a 
$U_v(\hat{\mathfrak sl}_r)\otimes 
U_v(\hat{\mathfrak sl}_m)$-module, we use the following coproducts. 
\[
\begin{array}{rl}
\Delta^{(l)}(f_{\bar i})\!\!&=f_{\bar i}\otimes 1 + 
v^{-h_{\bar i}}\otimes f_{\bar i}\\
&\\
\Delta^{(r)}(f_{\bar i})\!\!&=f_{\bar i}\otimes 1 + 
v^{h_{\bar i}}\otimes f_{\bar i}
\end{array}
\]
Note that $\Delta^{(l)}$ (resp. $\Delta^{(r)}$) is 
obtained from Lusztig's coproduct (resp. $\Delta'$) 
by reversing the order of the 
tensor product. The only reason we use them here is that we are 
more familiar with semi-infinite wedges which are infinite 
to the right. $\Delta^{(l)}$ behaves well for the bases at $v=\infty$. 
On the other hand, $\Delta^{(r)}$ behaves well for 
Kashiwara's lower crystal bases. Thus we are 
mostly concerned with $\Delta^{(l)}$ to work with 
${\cal F}_{v^{-1}}^{-\tilde\gamma}$. 

If $i_k$ are in descending order, 
they are called {\bf normally ordered} wedges. 
The straightnening laws 
are explained in \cite[(Ri)(Rii)(Riii)(Riv)]{U}. 
The normally ordered 
semi-infinite wedges of charge $c$ form a basis of 
${\cal F}_c$ \cite[Proposition 4.1]{TU}. 
For a normally ordered wedge, we locate them on 
an abacus with $rm$ runners. On each runner, larger numbers 
appear in upper location, and the row containing $1$ is read 
$1,\dots,rm$ from left to right. 
We divide the set of these 
runners into $m$ blocks. Then we have $m$ abacuses each of 
which has $r$ runners. By reading $i_k$'s in each block, 
we have $m$ semi-infinite wedges. We now assume that these are 
of the form 
$u_I^{(k)}:=u_{j_1^{(k)}}\wedge u_{j_2^{(k)}}\wedge\cdots$ 
such that 
$j_i^{(k)}=-\tilde\gamma_k\!-\!i\!+\!1$ for all $k$ and 
$i\!>\!>\!0$. 
We then identify $u_I^{(k)}$ with a multipartition 
$\lambda^{(k)}$ by 
$j_i^{(k)}=-\tilde\gamma_k\!+\!{\lambda^{(k)}}_i\!-\!i\!+\!1$. 
We identify ${\cal F}_{v^{-1}}^{-\tilde\gamma}$ with 
the subspace of 
${\cal F}_c$ ($c=-\sum \tilde\gamma_k$) 
spanned by the wedges $u_I$ whose 
$u_I^{(k)}$ have this form. 
This correspondence from normally ordered wedges to 
multipartitions is compatible with the action of 
$U_v(\hat{\mathfrak sl}_r)$. (If we consider the usual abacus with 
$r$ runners, it is compatible with $U_v(\hat{\mathfrak sl}_m)$-action.) 
More precisely, for each $n$, 
we take $\tilde\gamma$ such that 
$-\tilde\gamma_k<\!<-\tilde\gamma_{k+1}$ for all $k$. 
Then the action of $U_v(\hat{\mathfrak sl}_r)$ on 
the multipartitions of size less than $n$ coincides with the 
action given to ${\cal F}_{v^{-1}}^{-\gamma}$. This 
follows from the definition of the coproduct 
$\Delta^{(l)}$. 

We are now in a position to introduce a bar 
operation on the space of semi-infinite 
wedges as in \cite[3.1]{U}. The definition is identical to 
the definition of the bar operation on level one modules introduced in 
\cite[Proposition 3.1]{LT}. The welldefinedness for level one 
modules is given in \cite[5.1-5.9]{LT2}. The same proof works 
for the semi-infinite wedges considered here. We also have 
that $f_{\bar i}$ commutes with the bar operation, and 
that the bar operation preserves the size of 
multipartitions. 

We state the properties of the bar operation due to 
Uglov. 
For level one modules, these are stated in 
\cite[Theorem 3.2, Theorem 3.3]{LT}. (The proof is given in 
\cite[7.1-7.4]{LT2}.) 

\begin{thm}
\label{Uglov}
{\rm (\cite[Theorem 3.2, Theorem 3.3]{U})}

\begin{description}
\item[{\rm(1)}] 
The bar operation preserves ${\cal F}_{v^{-1}}^{-\tilde\gamma}$. 

\item[{\rm(2)}] 
$\overline{f_{\bar i}\underline\lambda}=f_{\bar i}
\overline{\underline\lambda}$, and $\overline{\underline\emptyset}
=\underline\emptyset$. In particular, 
the bar operation is an extension of the bar operation defined on 
${\cal M}_{v^{-1}}^{-\gamma}$. 

\item[{\rm(3)}] 
For each $n$, we take $\tilde\gamma$ as before. Then 
for multipartitions of size less than $n$, we have that 
$\overline{\underline\lambda}$ 
has the form 
$\underline\lambda+
\sum_{\underline\mu<\underline\lambda}
\alpha_{\underline\lambda,\underline\mu}(v)\underline\mu$. 
\end{description}
\end{thm}

The validity of $\overline{\underline\emptyset}
=\underline\emptyset$ comes from the facts that 
the bar operation preserves the subspace 
${\cal F}_{v^{-1}}^{-\tilde\gamma}$, and $\underline\emptyset$ 
is the unique multipartiton of the minimum size. 
The straightening laws show the unitriangularity of the 
bar operation. Note that the dominance order in \cite{U} corresponds 
to the reversed dominance order here. 
This triangularity also gives an algorithm 
to compute canonical basis on higher level modules. 
Thus it also computes decomposition numbers of cyclotomic Hecke algebras 
of type $G(m\!,1\!,n)$ over the field of complex numbers \cite{A1}. 

\bigskip
Since this theorem gives the required property of the canonical basis 
elements in question, we have reached the following theorem, 
which verifies the conjecture. 

\begin{thm} 
\label{main theorem}
$D^{\underline\lambda}\ne0$ if and only if 
$\underline\lambda$ is a Kleshchev multipartition. 
\end{thm}

\bigskip
\noindent
Tokyo University of Mercantile Marine,

\noindent
Etchujima 2-1-6, Koto-ku, Tokyo 135-8533, Japan 

\noindent
ariki@ipc.tosho-u.ac.jp


\begin{thebibliography}{aaa9}

\bibitem[A1]{A1}
S.Ariki, On the decomposition numbers of the Hecke algebra of 
$G(m,1,n)$, {\rm J.Math.Kyoto Univ.} {\bf 36} (1996), 789-808. 

\bibitem[A2]{A2}
S.Ariki, Representations over Quantum Algebras of type 
$A_{r-1}^{(1)}$ and Combinatorics of Young Tableux, 
{\rm Sophia University Lecture Notes Series (in Japanese)}, 
to appear. 

\bibitem[AM]{AM}
S.Ariki and A.Mathas, The number of simple modules of 
the Hecke algebras of type $G(r,1,n)$, {\rm Math.Zeit.}, 
to appear. 

\bibitem[BM]{BM}
M.Brou\'e and G.Malle, Zyklotomische Heckealgebren, 
{\rm Ast\'erisque} {\bf 212} (1993), 119-189. 

\bibitem[DJ1]{DJ1}
R.Dipper and G.James, Representations of Hecke algebras of 
general linear groups, 
{\rm Proc.London Math.(3)} {\bf 52} (1986), 20-52. 

\bibitem[DJ2]{DJ2}
R.Dipper and G.James, Representations of Hecke algebras of 
type $B_n$, {\rm Journal of Algebra} {\bf 146} (1992), 454-481. 

\bibitem[DJM]{DJM}
R.Dipper, G.James and A.Mathas, Cyclotomic $q$-Schur algebras, 
{\rm Math.Zeit.}, to appear. 

\bibitem[DJM']{DJM'}
R.Dipper, G.James and E.Murphy, Hecke algebras of type $B_n$ at 
roots of unity, {\rm Proc.London Math.Soc.(3)} {\bf 70} (1995), 
505-528. 

\bibitem[DM]{DM}
R.Dipper and A.Mathas, Morita equivalences of Ariki-Koike 
algebras, in preparation. 

\bibitem[GL]{GL}
J.J.Graham and G.I.Lehrer, Cellular algebras, 
{\rm Invent.Math.} {\bf 123} (1996), 1-34. 

\bibitem[G'L']{G'L'}
I.Grojnowski and G.Lusztig, A comparison of bases of quantized enveloping 
algebras, {\rm Contemp.Math.} {\bf 153} (1993), 11-19.

\bibitem[JMMO]{JMMO}
M.Jimbo, K.C.Misra, T.Miwa and M.Okado, Combinatorics of 
representations of $U_q(\hat{sl}(n)$ at $q=0$, 
{\rm Comm.Math.Phys.} {\bf 136} (1991), 543-566. 

\bibitem[Ka]{Ka}
M.Kashiwara, On crystal bases of the q-analogue of universal 
enveloping algebras, {\rm Duke Math.J.} {\bf 63} (1991), 465-516.

\bibitem[KMS]{KMS}
M.Kashiwara, T.Miwa and E.Stern, Decomposition of $q$-deformed 
Fock spaces, {\rm Selecta Math.} {\bf New Series 1} (1995), 787-805. 

\bibitem[Lamb]{Lamb}
S.Lambropoulou, Knot theory related to generalized and 
cyclotomic Hecke algebras of type B, {\rm Journal of knot theory and 
its ramifications}, to appear.

\bibitem[LT1]{LT}
B.Leclerc and J-Y.Thibon, Canonical bases of $q$-deformed Fock spaces, 
{\rm IMRN} {\bf 9} (1996), 447-455.

\bibitem[LT2]{LT2}
B.Leclerc and J-Y.Thibon, Littlewood-Richardson coefficients 
and Kazhdan-Lusztig polynomials, 
{\bf math.QA/9809122}. 

\bibitem[LLT]{LLT}
A.Lascoux, B.Leclerc and J-Y.Thibon, Hecke algebras at roots of 
unity and crystal bases of quantum affine algebras, 
{\rm Comm.Math.Phys.} {\bf 181} (1996), 205-263. 

\bibitem[L1]{L0}
G.Lusztig, Quivers, perverse sheaves, and quantized enveloping 
algebras, {\rm J.A.M.S.} {\bf 4} (1991), 365-421. 

\bibitem[L2]{L1}
G.Lusztig, Introduction to Quantum Groups, {\rm Progress in Math.} 
{\bf 110} (1993), Birkh\"auser.

\bibitem[L3]{L2}
G.Lusztig, Canonical basis and Hall algebras, 
{\rm Representation Theories and Algebraic Geometry, A.Broer and 
A.Daigneault eds.}, {\rm NATO ASI series C} {\bf 514} (1998), 365-399. 

\bibitem[MM]{MM}
T.Misra and K.C.Miwa, Crystal base for the basic representation 
of $U_q(\hat{sl}(n))$, 
{\rm Comm.Math.Phys.} {\bf 134} (1990), 79-88. 

\bibitem[N1]{N1}
H.Nakajima, Instantons on ALE spaces, quiver varieties, and 
Kac-Moody algebras, 
{\rm Duke Math.J.} {\bf 76} (1994), 365-416.

\bibitem[N2]{N2}
H.Nakajima, Quiver varieties and Kac-Moody algebras, 
{\rm Duke Math.J} {\bf 91} (1998), 515-560. 

\bibitem[TU]{TU}
K.Takemura and D.Uglov, Representations of the quantum toroidal 
algebra on highest weight modules of the quantum affine algebra 
of type ${\mathfrak gl}_N$, 
{\bf math.QA/9806134}. 

\bibitem[U]{U}
D.Uglov, Canonical bases of higher-level $q$-deformed Fock spaces, 
short version in {\bf math.QA/9901032}; full version in 
{\bf math.QA/9905196}. 

\bibitem[Vig1]{Vig1}
M-F.Vigneras, A propos d\'une conjecture de Langlands modulaire, 
{\rm Finite Reductive Groups, related structures and representations, 
M.Cabanes eds.} (1996), Birkh\"auser. 

\bibitem[Vig2]{Vig2}
M-F.Vigneras, Induced $R$-representations of $p$-adic reductive 
groups, {\rm Selecta Mathematica}, {\bf New Series 4} (1998), 549-623. 

\bibitem[Vig3]{Vig3}
M-F.Vigneras, private communication. 

\bibitem[VV]{VV}
M.Varagnolo and E.Vasserot, On the decomposition matrices of 
the quantized Schur algebra, {\rm Duke Math.J.}, to appear, 
{\bf math.QA/9803023}. 

\end{thebibliography}
\end{document}